\documentclass[12pt,a4paper]{article}
\usepackage{graphicx}
\usepackage[latin1]{inputenc}
\usepackage[T1]{fontenc}
\usepackage[english]{babel}
\usepackage{amssymb,amsmath}

\newcommand{\augmentelargeur}[1]{
\addtolength{\evensidemargin}{-#1}
\addtolength{\oddsidemargin}{-#1}
\addtolength{\textwidth}{#1}
\addtolength{\textwidth}{#1}
}
\augmentelargeur{5mm}

\newcommand{\thmheadercommand}[1]{\textbf{\scshape{}#1}}

\def\d{{\mathrm{d}}}
\def\N{{\mathbb{N}}}
\def\Z{{\mathbb{Z}}}

\renewcommand{\geq}{\geqslant}
\renewcommand{\leq}{\leqslant}

\def\eps{\varepsilon}
\renewcommand{\epsilon}{\varepsilon}
\renewcommand{\phi}{\varphi}

\newcommand{\abs}[1]{\left|\mskip1mu#1\right|}
\newcommand{\norm}[1]{\left\|#1\right\|}

\newcommand{\presgroup}[2]{\left\langle\,#1 \mid  #2\,\right\rangle}

\newcounter{prop}
\newcounter{defi}
\newcounter{thm}
\newcounter{lem}

\newenvironment{dem}[1][]{\noindent{\thmheadercommand{Proof#1}}\,\,--\,\,}{$\square$\medskip}

\newenvironment{enonce}[1]{\medskip\noindent{\thmheadercommand{#1}}\,\,--\,\,\begin{slshape}}{\end{slshape}\medskip}

\newenvironment{enonce2}[1]{\medskip\noindent{\thmheadercommand{#1}}\,\,--\,\,}{\medskip}

\newenvironment{defi}[1][]{\refstepcounter{prop}
\begin{enonce}{Definition \theprop{}#1}}{\end{enonce}}
\newenvironment{prop}[1][]{\refstepcounter{prop}
\begin{enonce}{Proposition \theprop{}#1}}{\end{enonce}}
\newenvironment{thm}[1][]{\refstepcounter{prop}
\begin{enonce}{Theorem \theprop{}#1}}{\end{enonce}}
\newenvironment{lem}[1][]{\refstepcounter{prop}
\begin{enonce}{Lemma \theprop{}#1}}{\end{enonce}}
\newenvironment{cor}[1][]{\refstepcounter{prop}
\begin{enonce}{Corollary \theprop{}#1}}{\end{enonce}}

\newenvironment{rem}[1][]{\refstepcounter{prop}
\begin{enonce2}{Remark \theprop{}#1}}{\end{enonce2}}

\newenvironment{defi*}[1][]{
\begin{enonce}{Definition#1}}{\end{enonce}}
\newenvironment{prop*}[1][]{
\begin{enonce}{Proposition#1}}{\end{enonce}}
\newenvironment{thm*}[1][]{
\begin{enonce}{Theorem#1}}{\end{enonce}}
\newenvironment{lem*}[1][]{
\begin{enonce}{Lemma#1}}{\end{enonce}}
\newenvironment{cor*}[1][]{
\begin{enonce}{Corollary#1}}{\end{enonce}}
\newenvironment{ex*}[1][]{
\begin{enonce}{Example#1}}{\end{enonce}}
\newenvironment{exo*}[1][]{
\begin{enonce2}{Exercise#1}}{\end{enonce2}}
\newenvironment{rem*}[1][]{
\begin{enonce2}{Remark#1}}{\end{enonce2}}

\title{Growth exponent of generic groups}

\author{Yann Ollivier}

\begin{document}
\maketitle

\def\d{\partial}

\begin{abstract}
In~\cite{GdlH}, Grigorchuk and de la Harpe ask if there are many groups
with growth exponent close to that of the free group with the same number
of generators. We prove that this is in fact the case for a generic
group (in the density model of random groups). Namely, for every positive
$\eps$, the property of having growth exponent at least $1-\eps$ (in base
$2m-1$ where $m$ is the number of generators) is generic in this model.
This extends in particular a theorem of Shukhov.

More generally, we prove that the growth exponent does not change much
through a random quotient of a torsion-free hyperbolic group.
\end{abstract}

The growth exponent is a very natural quantity associated to a group
presentation. It measures the rate of growth of the balls in the group
with respect to some given set of generators. Namely, let
$G=\presgroup{a_1,\ldots,a_m}{R}$ be a finitely generated group. For $L\geq 0$ let
$B_L\subset G$ be the set of elements of norm at most $L$ with respect to
this generating set. The
\emph{growth exponent} of $G$ (sometimes called \emph{entropy}) with respect to this set of generators is
\[
g=\lim_{L\rightarrow\infty} \frac1L \log_{2m-1} \abs{B_L}
\]

The maximal value of $g$ is achieved for the free group $F_m$, for which
$g=1$. The limit exists thanks to the submultiplicativity property
$\abs{B_{L+L'}}\leq \abs{B_L}\abs{B_{L'}}$. This implies in particular
that for any $L$ we have $\abs{B_L}\geq (2m-1)^{gL}$.

Growth exponent of groups, first introduced by Milnor, are related to
many other properties, for example in Riemannian geometry, dynamical
systems and of course combinatorial group theory. We refer
to~\cite{GdlH}, \cite{Har00} (chapters VI and VII), or \cite{Ver00} for
some surveys and applications, or to~\cite{Har02}, \cite{Wil04},
\cite{Kou98} for research more oriented towards uniform growth exponents.

\bigskip

The authors of~\cite{GdlH} ask if there are lots of families of groups
whose growth exponents gets arbitrarily close to the maximal value $1$.
An example of such groups is given in~\cite{Shu}: Shukhov proves therein
that if a group presentation satisfies the $C'(1/6)$ small cancellation
condition and has ``not too many'' relators (in a precise sense), then if
the relators are long the growth exponent of the group so presented is
very close to $1$.

We prove that in fact, lots of groups share this property: having growth
exponent at least $1-\eps$ is a \emph{generic} property in the density
model of random groups (introduced by Gromov in~\cite{Gro93}, and which
we recall precisely in section~\ref{randomgroups} below). The density
model allows a precise control of the quantity of relations put in the
random group, as examplified by the phase transition proven
in~\cite{Gro93}: below density $1/2$, random groups are very probably
infinite and hyperbolic, and, above density $1/2$, very probably trivial.
In this framework, our theorem reads:

\begin{thm}
\label{main}
Let $d<1/2$ be a density parameter and let $G$ be a random group on
$m\geq 2$ generators at density $d$ and length $\ell$.

Then, for any $\eps>0$, the probability that the growth exponent of $G$
is at least $1-\eps$ tends to $1$ as $\ell\rightarrow\infty$.
\end{thm}

When $d<1/12$ this is a consequence of Shukhov's theorem: indeed for
densities at most $1/12$, random groups satisfy the $C'(1/6)$ small
cancellation condition. But for larger densities they do not any more, and
so the theorem really provides a large class of new groups with large
growth exponent.

Random groups at length $\ell$ look like free groups at scales lower than
$\ell$, and so the cardinality of balls of course grows with exponent $1$
at the beginning. However growth is an asymptotic invariant, and the
geometry of random groups at scale $\ell$ is highly non-trivial, so the
theorem cannot be interpreted by simply saying that random groups look
like free groups at small scales.

\bigskip

More generally, we show that for torsion-free hyperbolic groups, the
growth exponent is stable in the following sense: if we randomly pick
elements in the group and quotient by the normal subgroup they generate
(the so-called \emph{geodesic} model of random quotient as opposed to
randomly picking words in the generators; see details below), then the
growth exponent stays almost unchanged (unless we killed too many
elements and get the trivial group).

The study of random quotients of hyperbolic groups arises
naturally from the knowledge that a random group (a random quotient of
the free group) is hyperbolic: one can wonder whether a random
quotient of a hyperbolic group stays hyperbolic. The answer
from~\cite{Oll03} is yes (see section~\ref{randomgroups} below for
details) up to some critical density equal to $g/2$ where $g$ is the
growth exponent of the initial group; above this critical density the
random quotient collapses. In this framework our second theorem reads:

\begin{thm}
\label{genmain}
Let $G_0$ be a non-elementary torsion-free hyperbolic group of growth
exponent $g$. Let $d<g/2$. Let $G$ be a random quotient of $G_0$ by
\emph{geodesic} words at length $\ell$.

Then, for any $\eps>0$, with probability tending to $1$ as
$\ell\rightarrow\infty$, the growth exponent of $G$ lies between $g-\eps$
and $g$.
\end{thm}

Of course, Theorem~\ref{main} is just Theorem~\ref{genmain} applied to a free
group.

\begin{rem}
The proof of Theorem~\ref{genmain} only uses the two following facts: that
the random quotient axioms of~\cite{Oll03} are satisfied, and that there is
a local-to-global principle for growth in the random quotient. So in
particular the result holds under slightly weaker conditions than
torsion-freeness of $G_0$, as described in~\cite{Oll03} (``harmless
torsion'').
\end{rem}

\paragraph{Locality of growth in hyperbolic groups.} As one of our tools
we use a result about locality of growth in hyperbolic groups (see the
Appendix).  Growth is an asymptotic invariant, and large relations in a
group can change it noticeably. But in hyperbolic groups, if the
hyperbolicity constant is known, it is only necessary to evaluate growth
in some ball in the group to get a bound for growth of the group (see
Proposition~\ref{localgr} in the Appendix). So in particular, growth of a
given hyperbolic presentation is computable.

In the case of random quotients by relators of length $\ell$, this
principle shows that it is necessary to check growth up to words of
length at most $A\ell$ for some large constant $A$ (which depends on
density and actually tends to infinity when $d$ is close to the critical
density), so that geometry of the quotient matters up to scale $\ell$
(including the non-trivial geometry of the random quotient at this scale)
but not at higher scales.

This result may have independent interest.

\paragraph{About the proofs, and about cogrowth.} The proofs presented
here make heavy use of the terminology and results from~\cite{Oll03}. We
have included a reminder (section~\ref{reminder}) so that this paper is
self-contained.

This paper comes along with a ``twin'' paper about \emph{co}growth of random
groups (\cite{Oll04}). Let us insist that, although the inspiration for
these two papers is somewhat the same (use some locality principle and
count van Kampen diagrams), they mostly differ in detail, except for the
reminder from~\cite{Oll03} which is identical. Especially, the proof of
the locality principle for growth and cogrowth is not at all the same.
The counting of van Kampen diagrams begins similarly but soon diverges as
we are not evaluating the same things eventually. And we do not work in
the same variant of the density model: for growth we use the geodesic
variant, whereas for cogrowth we use the word variant (happily these two
variants coincide in he case of a free group, that is, for ``plain''
random groups).

\paragraph{Acknowledgments.} I would like to thank Étienne Ghys and
Pierre Pansu for helpful discussions and many comments on the text.
Lots of the ideas presented here emerged during my stay at the École
normale supérieure de Lyon in Spring 2003, at the invitation of Damien
Gaboriau and Étienne Ghys. I am very grateful to all the team of the
math department there for their warmth at receiving me.

\section{Definitions and notations}

\subsection{Random groups and density}
\label{randomgroups}

The interest of random groups is twofold: first, to study which
properties of groups are \emph{generic}, i.e.\ shared by a large
proportion of groups; second, to provide examples of new groups with
given properties. This article falls under both approaches.

A random group is given by a random presentation, that is, the quotient
of a free group $F_m=\langle a_1,\ldots,a_m\rangle$ by (the normal
closure of) a randomly chosen set $R\subset F_m$.
Defining a random group is giving a law for the random set $R$.

More generally, a random quotient of a group $G_0$ is the quotient of
$G_0$ by (the normal closure of) a randomly chosen subset $R\subset G_0$.

The philosophy of random gorups was introduced by Gromov in~\cite{Gro87}
through a statement that ``almost every group is hyperbolic'', the proof of
which was later given by Ol'shanski\u{\i} (\cite{Ols92}) and independently
by Champetier (\cite{Ch95}). Gromov later defined the density model
in~\cite{Gro93}, in order to precisely control the quantity of relators
put in a random group.

Since then random groups have gained broad interest and are connected
to lots of topics in geometric or combinatorial group theory (such
as the isomorphism problem, property T, Haagerup property, small
cancellation, spectral gaps, the Baum-Connes conjecture...), especially
since Gromov used them (\cite{Gro03}) to build a counter-example to the
Baum-Connes conjecture with coefficients (see also~\cite{HLS}). We refer
to~\cite{Gh} for a general discussion on random groups.

\bigskip

We now define the density model of random groups. In this model the
random set of relations $R$ depends on a density parameter $d$: the
larger $d$, the larger $R$. This model exhibits a phase transition
between infiniteness and triviality depending on the value of $d$;
moreover, in the infinite phase some properties of the resulting group
(such as the rank, property T or the Haagerup property) do differ
depending on $d$, hence the interest of this model.

\begin{defi}[ (Density model of random groups or quotients, geodesic
variant)]
Let $G_0$ be a group generated by the elements $a_1^{\pm
1},\ldots,a_m^{\pm 1}$ ($m\geq 2$). Let $B_\ell\subset G_0$ be the ball
of radius $\ell$ in $G_0$ with respect to this generating set.

Let $0\leq d\leq 1$ be a density parameter.

Let $R$ be a set of $(2m-1)^{d\ell}$ randomly chosen elements of
$B_\ell$, uniformly and independently picked in $B_\ell$.

We call the group $G=G_0/\langle R\rangle$ a \emph{random quotient} of
$G_0$ by geodesic words, at density $d$ and length $\ell$.

In case $G_0$ is the free group $F_m$ we simply call $G$ a \emph{random
group}.
\end{defi}

In this definition, we can also replace $B_\ell$ by the sphere $S_\ell$
of elements of norm exactly $\ell$, or by the annulus of elements of norm
between $\ell$ and $\ell+C$ for some constant $C$: this does not affect
our theorems.

Another variant (the word variant) of random groups consists in taking
for $R$ a set of reduced (or plain) words in the generators $a_i^{\pm
1}$, which leads to a different probability distribution. Fortunately in
the case of the free group, there is no difference between taking at
random elements in $B_\ell$ or reduced words, so that the notions of
random group and of a generic property of groups are well-defined anyway.
Quotienting by elements rather than words seems better suited to control
the growth of the quotient.

The interest of the density model was established by the following
theorem of Gromov, which shows a sharp phase transition between infinity
and triviality of random groups.

\begin{thm}[ (M.~Gromov, \cite{Gro93})]
Let $d<1/2$. Then with probability tending to $1$ as $\ell$ tends to
infinity, random groups at density $d$ are infinite hyperbolic.

Let $d>1/2$. Then with probability tending to $1$ as $\ell$ tends to
infinity, random groups at density $d$ are either $\{e\}$ or $\Z/2\Z$.
\end{thm}

(The occurrence of $\Z/2\Z$ is of course due to the case when $\ell$ is
even and we we take elements in the sphere $S_\ell$; this disappears if
one takes elements in $B_\ell$, or of length between $\ell$ and $\ell+C$
with $C\geq 1$.)

Basically, $d\ell$ is to be interpreted as the ``dimension'' of the
random set $R$ (see the discussion in~\cite{Gro93}). As an illustration,
if $L< 2d\ell$ then very probably there will be two relators in $R$
sharing a common subword of length $L$. Indeed, the dimension of the
couples of relators in $R$ is $2d\ell$, whereas sharing a common subword
of length $L$ amounts to $L$ ``equations'', so the dimension of those
couples sharing a subword is $2d\ell-L$, which is positive if $L<2d\ell$.
This ``shows'' in particular that at density $d$, the small cancellation
condition $C'(2d)$ is satisfied.

\bigskip

Since a random quotient of a free group is hyperbolic, one can wonder if
a random quotient of a hyperbolic group is still hyperbolic. The answer
is basically yes, and for the geodesic variant, the critical density in
this case is linked to the growth exponent of the initial group.

\begin{thm}[ (Y.~Ollivier, \cite{Oll03})]
\label{thmrq}
Let $G_0$ be a non-elementary,
torsion-free hyperbolic group, generated by the elements $a_1^{\pm
1},\ldots,a_m^{\pm 1}$, with cogrowth exponent $g$.
Let $0\leq d \leq 1$ be a density parameter.

If $d<g/2$, then a random quotient of $G_0$ by random geodesic
words at density $d$ is infinite hyperbolic, with probability tending to
$1$ as $\ell$ tends to infinity.

If $d>g/2$, then a random quotient of $G_0$ by random geodesic
words at density $d$ is either $\{e\}$ or $\Z/2\Z$, with probability tending to
$1$ as $\ell$ tends to infinity.
\end{thm}

This is the context in which Theorem~\ref{genmain} is to be understood.

\subsection{Hyperbolic groups and isoperimetry of van Kampen diagrams}
\label{defiso}

Let $G$ be a group given by the finite presentation
$\presgroup{a_1,\ldots,a_m}{R}$. Let $w$ be a word in the $a_i^{\pm
1}$'s. We denote by $\abs{w}$ the number of letters of $w$, and by
$\norm{w}$ the distance from $e$ to $w$ in the Cayley graph of the
presentation, that is, the minimal length of a word representing the same
element of $G$ as $w$.

Let $\lambda$ be the maximal length of a relation in $R$.

\def\d{\partial}
\def\A{\mathcal{A}}

We refer to~\cite{LS} for the definition and basic properties of van
Kampen diagrams. Remember that a word represents the neutral element of
$G$ if and only if it is the boundary word of some van Kampen diagram. If
$D$ is a van Kampen diagram, we denote its number of faces by $\abs{D}$
and its boundary length by $\abs{\d D}$.

It is known (see for example~\cite{Sho}) that $G$ is hyperbolic if and only if there
exists a constant $C_1>0$ such that for any (reduced) word $w$
representing the neutral element of $G$, there exists a van
Kampen diagram with boundary word $w$, and with at most $\abs{w}/C_1$ faces.
This can be reformulated as: for any word $w$ representing the neutral
element of $G$, there exists a van Kampen diagram with boundary word $w$
satisfying the isoperimetric inequality
\[
\abs{\d D}\geq C_1 \abs{D}
\]

We are going to use a \emph{homogeneous} way to write this inequality.
The above form compares the boundary length of a van Kampen diagram to
its number of faces. This amounts to comparing a length with a number,
which is not very well-suited for geometric arguments, especially when
dealing with groups having relations of very different lengths.

So let $D$ be a van Kampen diagram w.r.t.\ the presentation and define the
\emph{area} of $D$ to be
\[
\A(D)=\sum_{f\text{ face of }D} \abs{\d f}
\]
which is also the number of external edges (not couting ``filaments'') plus twice the number of
internal ones. This has, heuristically speaking, the homogeneity of a
length.

It is immediate to see that if $D$ satisfies $\abs{\d D}\geq C_1
\abs{D}$, then we have $\abs{\d D}\geq C_1\,\A(D)/\lambda$ (recall
$\lambda$ is the maximal length of a relation in the presentation).
Conversely, if $\abs{\d D}\geq C_2\,\A(D)$, then $\abs{\d D}\geq
C_2\abs{D}$. So we can express the isoperimetric inequality using $\A(D)$
instead of $\abs{D}$.

Say a diagram is \emph{minimal} if it has minimal area for a given
boundary word. So $G$ is hyperbolic if and only if there exists a
constant $C>0$ such that every minimal van Kampen diagram satisfies the
isoperimetric inequality
\[
\abs{\d D}\geq C\,\A(D)
\]
and of course we necessarily have $C\leq 1$ (since the edges making up
$\d D$ are taken into account in $\A(D)$).

This formulation is homogeneous, that is, it compares a length to a
length. This inequality is the one that naturally arises in $C'(\alpha)$
small cancellation theory (with $C=1-6\alpha$) as well as in random
groups at density $d$ (with $C=\frac12-d$). So in these contexts the value of
$C$ is naturally linked with some parameters of the presentation.

This kind of isoperimetric inequality is also the one appearing in the
assumptions of Champetier in~\cite{Ch93}, in random quotients of
hyperbolic groups (cf.~\cite{Oll03}) and in the (infinitely presented)
limit groups constructed by Gromov in~\cite{Gro03}.  So we think this is
the right way to write the isoperimetric inequality when the lengths of
the relators are very different.

\bigskip

This formulation has yet another advantage in that the relationship
between the isoperimetric constant $C$ and the hyperbolicity constant
$\delta$ (Rips' constant for thinness of triangles, see~\cite{GH}
or~\cite{Sho} for definitions) is more elegant. Indeed, we have the
following:

\begin{prop}
\label{isomdelta}
Suppose that the hyperbolic group $G$ given by some finite presentation
satisfies the isoperimetric inequality
\[
\abs{\d D}\geq C\,\A(D)
\]
for all minimal van Kampen diagrams $D$, for some constant $C>0$.

Let $\lambda$ be the maximal length of a relation in the
presentation. Then the hyperbolicity constant $\delta$ of $G$ satisfies
\[
\delta\leq 12\lambda/C^2
\]
\end{prop}

The point to note is of course that both $\delta$ and $\lambda$ are
lengths whereas $C$ is purely numerical (between $0$ and $1$), hence
homogeneity of the result.  This will be crucial to deal with random
quotients of hyperbolic groups, in which the relations are of very
different lengths.

\medskip

\begin{dem}
Actually the proof of this is strictly included in~\cite{Sho}
(Theorem~2.5). Indeed, what the authors of~\cite{Sho} prove is always of the form ``the
number of edges in $D$ is at least something, so the number of faces of
$D$ is at most this thing divided by $\rho$'' (in their notation $\rho$
is the maximal length of a relation). Reasoning directly with
the number of edges instead of the number of faces $\abs{D}$
simplifies their arguments. But $\A(D)$ is simply twice the number of
internal edges of $D$ plus the number of boundary edges of $D$, so it is
greater than the number of edges of $D$.

So (keeping their notations) simply by removing the seventh sentence
their Lemma~2.6 (evaluating the number of $2$-cells by dividing the
number of $2$-cells by $\rho$) gives a new Lemma~2.6 which reads (we
stick to their notations in the framework of their proving Theorem~2.5)

\begin{lem*}[ 2.6 of~\cite{Sho}]
If $\eps>\rho$, then there is a constant $C_1$ depending solely on $\eps$,
such that the number of $1$-cells in $N(\theta)$ is at least
$\ell(\theta)\eps/\rho-C_1$. Namely we can set $C_1=\eps(\eps+\rho)/\rho$.
\end{lem*}

Similarly, removing the last sentence of their proof of Lemma~2.7 we get
a new version of it:

\begin{lem*}[ 2.7 of~\cite{Sho}]
If $\eps>\rho$, there is a constant $C_2$ depending solely on $\eps$ such
that 
\[
\A(D)> (\alpha+\beta+\gamma)\eps/\rho-C_2+2r/\rho
\]
where $\A(D)$ is the area of the diagram $D$. Namely we can set
$C_2=3C_1+4\eps+2$.
\end{lem*}

We insist that those modified lemmas are obtained by \emph{removing} some
sentences in their proofs, and that there really is nothing to modify.

We still have to re-write the conclusion. In their notation $\alpha$,
$\beta$ and $\gamma$ are (up to $4\eps$) the lengths of the sides of some
triangle which, by contradiction, is supposed not to be $r$-thin (we want
to show that if $r$ is large enough, then every triangle is $r$-thin).

The assumption $\abs{\d D}\geq C\,\A(D)$ reads
\[
\A(D)\leq (\alpha+\beta+\gamma)/C+12\eps/C
\]
Combining this inequality and the result of Lemma~2.7, we have
\[
(\alpha+\beta+\gamma)\eps/\rho-C_2+2r\leq
(\alpha+\beta+\gamma)/C+12\eps/C
\]
Now set $\eps=\rho/C$. We thus obtain
\[
2r\leq 12\rho/C^2+C_2
\]
where we recall that
$C_2=3C_1+4\eps+2=3\eps(\eps+\rho)/\rho+4\eps+2=\rho(3/C^2+7/C)+2$ with
our choice of $\eps$. Since $\rho\geq 1$ (unless $G$ is free in which
case there is nothing to prove) and necessarily $C\leq 1$ we have
$7/C\leq 7/C^2$ and $2\leq 2\rho/C^2$ and so finally
\[
2r\leq 12\rho/C^2+12\rho/C^2
\]
hence the conclusion, remembering that our $\delta$ and $\lambda$
are~\cite{Sho}'s $r$ and $\rho$ respectively.
\end{dem}

\section{Growth of random quotients}
\label{mainsection}

We now turn to the main point of this paper, namely, evaluation of the
growth exponent of a random quotient of a group.

\subsection{Framework of the argument}

So let $G_0$ be a non-elementary torsion-free hyperbolic group given by the
finite presentation $G_0=\presgroup{a_1,\ldots,a_m}{Q}$. Let $g>0$ be the
growth exponent of $G_0$ with respect to this generating set. Let $B_\ell$
be the set of elements of norm at most $\ell$. Let $\lambda$ be the
maximal length of a relation in $Q$.

Let also $R$ be a randomly chosen set of $(2m-1)^{d\ell}$ elements of the
ball $B_\ell\subset G_0$, in accordance with the model of random
quotients we retained. Set $G=G_0/\langle R\rangle$, the random quotient
we are interested in. We will call the relators in $R$ ``new relators''
and those in $Q$ ``old relators''.

In the sequel, the phrase ``with overwhelming probability'' will mean
``with probability exponentially tending to $1$ as
$\ell\rightarrow\infty$ (depending on everything)''.

\bigskip

\def\B{\mathcal{B}}

Fix some $\eps>0$. We want to show that the growth exponent of $G$
is at least $g(1-\eps)$, with overwhelming probability.

We can suppose that the length $\ell$ is taken large enough so that, for
$L\geq \ell$, we have $(2m-1)^{gL}\leq \abs{B_L}\leq
(2m-1)^{g(1+\eps)L}$.

Let $\B_L$ be the ball of radius $L$ in $G$. We trivially have $\abs{\B_L}\leq\abs{B_L}$.

We will prove a lower bound for the cardinality of $\B_L$ for some well
chosen $L$, and then use Proposition~\ref{localgr}. In order to apply
this proposition, we first need an estimate of the hyperbolicity
constant of $G$.

\begin{prop}
With overwhelming probability, minimal van
Kampen diagrams of $G$ satisfy the isoperimetric
inequality
\[
\abs{\d D}\geq C\,\A(D)
\]
where $C>0$ is a constant depending on $G_0$ and the density
$d$ but not on $\ell$. In particular, the
hyperbolicity constant $\delta$ of $G$ is at most
$12\ell/C^2$. 
\end{prop}

\begin{dem}
This is a rephrasing of Proposition~32 of~\cite{Oll03}: With overwhelming probability, minimal van Kampen diagrams $D$ of the
random quotient $G$ satisfy the isoperimetric inequality
\[
\abs{\d D}\geq \alpha_1\ell\abs{D''}+\alpha_2\abs{D'}
\]
where $\alpha_1,\alpha_2$ are positive constants depending on $G_0$ and the
density parameter $d$ (but not on $\ell$), and $\abs{D''}$, $\abs{D'}$
are respectively the number of faces of $D$ bearing new relators (from
$R$) and old relators (from $Q$). Since new relators have length at most $\ell$
and old relators have length at most $\lambda$, by definition we have
$\A(D)\leq \ell\abs{D''}+\lambda\abs{D'}$
and so setting
$C=\min(\alpha_1,\alpha_2/\lambda)$ yields
\[
\abs{\d D}\geq C\,\A(D)
\]

The estimate of the hyperbolicity constant follows by
Proposition~\ref{isomdelta}.
\end{dem}

In particular, in order to apply Proposition~\ref{localgr} it is
necessary to control the cardinality of balls of radius roughly
$\ell/C^2+1/g$. More precisely, let $A\geq 500$ be such that $40/A\leq
\eps/2$. Set $L_0=24\ell/C^2+4/g$ and $L=AL_0$. We already trivially know that
$\abs{\B_{L_0}}\leq (2m-1)^{g(1+\eps)L_0}$. We will now show that, with
overwhelming probability, we have $\abs{\B_L}\geq (2m-1)^{g(1-\eps/2)L}$.
Once this is done we can conclude by Proposition~\ref{localgr}.

\bigskip

The strategy to evaluate the growth of the quotient $G$ of $G_0$ will be
the following: There are at least $(2m-1)^{gL}$ elements in $B_L$.  Some
of these elements are identified in $G$. Let $N$ be the number of
equalities of the form $x=y$, for $x,y\in B_L$, which hold in $G$ but did
not hold in $G_0$.  Each such equality decreases the number of elements
of $\B_L$ by at most $1$. Hence, the number of elements of norm $L$ in
$G$ is at least $(2m-1)^{gL}-N$. So if we can show for example that
$N\leq \frac12 (2m-1)^{gL}$, we will have a lower bound for the size of
balls in $G$.

So we now turn to counting the number of equalities $x=y$ holding in
$G$ but not in $G_0$, with $x,y\in B_L$. Each such
equality defines a (minimal) van Kampen diagram with boundary word
$xy^{-1}$, of boundary length at most $2L$. We will need the
properties of van Kampen diagrams of $G$ proven
in~\cite{Oll03}.

So, for the $\eps$ and $A$ fixed above, let $A'=2L/\ell$ and let $D$ be a
minimal van Kampen diagram of $G$, of boundary length at most
$A'\ell$. By the isoperimetric inequality $\abs{\d D}\geq C\A(D)$, we know
that the number $\abs{D''}$ of faces of $D$ bearing a new relator of $R$
is at most $A'/C$. So for all the sequel set
\[
K=A'/C
\]
which is the maximal number of new relators in the diagrams we have to
consider (which will also have area at most $K\ell$). Most importantly,
this $K$ does not depend on $\ell$.

\subsection{Reminder from~\cite{Oll03}}
\label{reminder}

In this context, it is proven in~\cite{Oll03} that the van Kampen diagram
$D$ can be seen as a ``van Kampen diagram at scale $\ell$ with respect to
the new relators, with equalities modulo $G_0$''. More precisely, this
can be stated as follows: (we refer to~\cite{Oll03} for the definition of
``strongly reduced'' diagrams; the only thing to know here is that for
any word equal to $e$ in $G$, there exists a strongly reduced van Kampen
diagram with this word as its boundary word).

\begin{prop}[ (\cite{Oll03}, section~6.6)]
\label{propdavKd}
Let $G_0=\presgroup{S}{Q}$ be a non-elementary hyperbolic group, let $R$ be a set of
words of length $\ell$, and consider the group $G=G_0/\langle R
\rangle=\presgroup{S}{Q \cup R}$.

Let $K\geq 1$ be an arbitrarily large integer and let $\eps_1,\eps_2>0$
be arbitrarily small numbers. Take $\ell$ large enough depending on $G_0,
K,\eps_1,\eps_2$.

Let $D$ be a van Kampen diagram with respect to the presentation
$\presgroup{S}{Q \cup R}$, which is strongly reduced, of area at most
$K\ell$. Let also $D'$ be the subdiagram of $D$ which is the union of the
$1$-skeleton of $D$ and of those faces of $D$ bearing relators in $Q$
(so $D'$ is a possibly non-simply connected van Kampen
diagram with respect to $G_0$), and suppose that
$D'$ is minimal.

We will call \emph{worth-considering} such a van Kampen diagram.

Let $w_1,\ldots,w_p$ be the boundary (cyclic) words of $D'$, so that each
$w_i$ is either the boundary word of $D$ or a relator in $R$.

\smallskip

Then there exists an integer $k\leq 3K/\eps_2$ and words
$x_2,\ldots,x_{2k+1}$ such that:
\begin{itemize}
\item Each $x_i$ is a subword of some cyclic word $w_j$;

\item As subwords of the $w_j$'s, the $x_i$'s are disjoint and their
union exhausts a proportion at least $1-\eps_1$ of the total length of
the $w_j$'s.

\item For each $i\leq k$, there exists words $\delta_1,\delta_2$ of
length at most $\eps_2 (\abs{x_{2i}}+\abs{x_{2i+1}})$ such that
$x_{2i}\delta_1 x_{2i+1} \delta_2=e$ in $G_0$.

\item If two words $x_{2i}$, $x_{2i+1}$ are subwords of the boundary words of
two faces of $D$ bearing the same relator $r^{\pm 1}\in R$, then, as
subwords of $r$, $x_{2i}$ and $x_{2i+1}$ are either disjoint or equal with
opposite orientations (so that the above equality reads $x\delta_1
x^{-1}\delta_2=e$).
\end{itemize}

The couples $(x_{2i},x_{2i+1})$ are called \emph{translators}.
Translators are called \emph{internal}, \emph{internal-boundary} or
\emph{boundary-boundary} according to whether $x_{2i}$ and $x_{2i+1}$
is a subword of some $w_j$ which is a relator in $R$ or the boundary word of
$D$.
\end{prop}

(There are slight differences between the presentation here and that
in~\cite{Oll03}. Therein, boundary-boundary translators did not have to be
considered: they were eliminated earlier in the process, before
section~6.6, because they have a positive contribution to boundary
length, hence always improve isoperimetry and do not deserve
consideration in order to prove hyperbolicity. Moreover, in~\cite{Oll03} we
further distinguished ``commutation translators'' for the kind of
internal translator with $x_{2i}=x_{2i+1}^{-1}$, which we need not do
here.)

Translators appear as dark strips on the following figure:

\begin{center}
\includegraphics{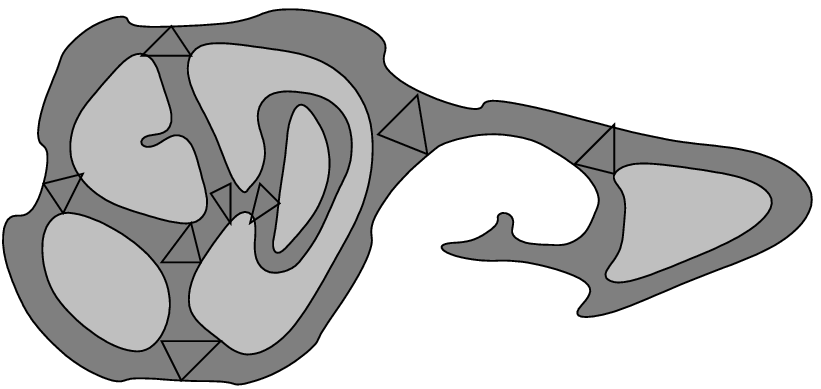}
\end{center}

\begin{rem}
\label{nosmalltrans}
Since there are at most $3K/\eps_2$ translators, the total length of the
translators $(x_{2i},x_{2i+1})$ for which $\abs{x_{2i}}+\abs{x_{2i+1}}\leq
\eps_3\ell$ is at most $3K\ell\eps_3/\eps_2$, which makes a proportion at
most $3\eps_3/\eps_2$ of the total length. So, setting
$\eps_3=\eps_1\eps_2/3$ and replacing $\eps_1$ with $\eps_1/2$, we can suppose that the union of the translators
exhausts a proportion at least $1-\eps_1$ of the total length of the
diagram, and that each translator $(x_{2i},x_{2i+1})$ satisfies
$\abs{x_{2i}}+\abs{x_{2i+1}}\geq\eps_1\eps_2\ell/6$.
\end{rem}

\begin{rem}
\label{numberofdavKds}
The number of ways to partition the words
$w_i$ into translators is at most $(2K\ell)^{12K/\eps_2}$, because each
$w_i$ can be determined by its starting- and endpoint, which can be given
as numbers between $1$ and $2K\ell$ which is an upper bound for the
cumulated length of the $w_i$'s (since the area of $D$ is at most
$K\ell$). For fixed $K$ and $\eps_2$ this grows subexponentially in
$\ell$.
\end{rem}

\begin{rem}
\label{orderboundary}
Knowing the words $x_i$, the number of possibilities for the boundary
word of the diagram is at most $(6K/\eps_2)!$ (choose which subwords
$x_i$ make the boundary word of the diagram, in which order), which does
not depend on $\ell$ for fixed $K$ and $\eps_2$.
\end{rem}

\bigskip

We need another notion from~\cite{Oll03}, namely, that of \emph{apparent
length} of an element in $G_0$. Apparent length is defined
in~\cite{Oll03} in a more general setting, with respect to a family of
measures on the group depending on the precise model of random quotient
at play. Here these are simply the uniform measures on the balls
$B_\ell$.  So we only give here what the definition amounts to in our
context. In fact we will not use here the full strength of this notion,
but we still need to define it in order to state results
from~\cite{Oll03}. 

Recall that in the geodesic model
of random quotients, the axioms of~\cite{Oll03} are satisfied with
$\beta=g/2$ and $\kappa_2=1$, by Proposition~20 of~\cite{Oll03}.

\def\L{\mathbb{L}}

\begin{defi}[ (Definition~36 of~\cite{Oll03})]
Let $x\in G_0$. Let $\eps_2>0$. Let $L$ be an integer. Let $p_L(xuyv=e)$
be the probability that, for a random element $y\in B_L$, there
exist elements $u,v\in G_0$ of norm at most $\eps_2(\norm{x}+L)$ such that
$xuyv=e$ in $G_0$.

The \emph{apparent length} of $x$ at test-length $L$ is
\[
\L_L(x)=-\frac{2}{g} \log_{2m-1} p_L(xuyv=e) - L
\]

The \emph{apparent length} of $x$ is
\[
\L(x)=\min\left(\norm{x}\,,\;\min_{0 \leq L
\leq K\ell} \L_L(x)\right)
\]
where we recall $\ell$ is the length of the relators in a random
presentation.
\end{defi}
%
%
%

%

\bigskip

We further need the notion of a \emph{decorated abstract van Kampen
diagram} (which was implicitly present in the free case when we mentioned
the probability that some diagram ``is fulfilled by random relators''),
which is inspired by Proposition~\ref{propdavKd}: it carries the
combinatorial information about how the relators and boundary word of a
diagram were cut into subwords in order to make the translators.

\def\D{\mathcal{D}}

\begin{defi}[ (Decorated abstract van Kampen diagram)]
\label{defdavKds}
Let $K\geq 1$ be an arbitrarily large integer and let $\eps_1,\eps_2>0$
be arbitrarily small numbers. Let $I_\ell$ be the cyclically ordered set of
$\ell$ elements.

A \emph{decorated abstract van Kampen diagram} $\D$ is the following data:
\begin{itemize}
\item An integer $\abs{\D}\leq K$ called its \emph{number of faces}.
\item An integer $\abs{\d \D}\leq K\ell$ called its \emph{boundary length}.
\item An integer $n\leq \abs{\D}$ called its \emph{number of distinct
relators}.
\item An application $r^\D$ from $\{1,\ldots,\abs{\D}\}$ to
$\{1,\ldots,n\}$; if $r^\D(i)=r^\D(j)$ we will say that \emph{faces $i$
and $j$ bear the same relator}.
\item An integer $k\leq 3K/\eps_2$ called the \emph{number of
translators} of $\D$.
\item For each integer $2\leq i\leq 2k+1$, a set of the form $\{j_i\}\times
I'_i$ where either $j_i$ is an integer between $1$ and $\abs{\D}$ and
$I'_i$ is an oriented cyclic subinterval of
$I_\ell$, or $j_i=\abs{\D}+1$ and $I'_i$ is a subinterval of $I_{\abs{\d
\D}}$;
this is called an
\emph{(internal) subword of the $j_i$-th face} in the first case, or a
\emph{boundary subword} in the second case.
\item For each integer $1\leq i \leq k$ such that $j_{2i}\leq \abs{\D}$, an integer between $0$ and
$4\ell$ called the \emph{apparent length of the $2i$-th subword}.
\end{itemize}
such that
\begin{itemize}
\item The sets $\{j_i\}\times I'_i$ are all disjoint and the cardinal
of their union is at least $(1-\eps_1)\left(\abs{\D}\ell+\abs{\d \D}\right)$.
\item For all $1\leq i\leq k$ we have $j_{2i}\leq j_{2i+1}$ (this can be
ensured by maybe swapping them).
\item If two faces $j_{2i}$ and $j_{2i+1}$ bear the same relator, then
either $I'_{2i}$ and $I'_{2i+1}$ are disjoint or are equal with opposite
orientations.
\end{itemize}
\end{defi}

This way, Proposition~\ref{propdavKd} ensures that any
worth-considering van Kampen diagram $D$ with respect to $G_0/\langle R
\rangle$ defines a decorated abstract van Kampen diagram $\D$ in the way
suggested by terminology (up to rounding the apparent lengths to the
nearest integer; we neglect this problem). We will say that $\D$ is
\emph{associated to} $D$.  Remark~\ref{numberofdavKds} tells that the
number of decorated abstract van Kampen diagrams grows subexponentially
with $\ell$ (for fixed $K$).

Given a decorated abstract van Kampen diagram $\D$ and $n$ given relators
$r_1,\ldots,r_n$, we say that these relators \emph{fulfill $\D$} if there
exists a worth-considering van Kampen diagram $D$ with respect to
$G_0/\langle r_1,\ldots,r_n\rangle$, such that the associated decorated
abstract van Kampen diagram is $\D$. Intuitively speaking, the relators
$r_1,\ldots,r_n$ can be ``glued modulo $G_0$ in the way described by
$\D$''.

So we want to study which diagrams can probably be fulfilled by random
relators in $R$. The main conclusion from~\cite{Oll03} is that these are
those with large boundary length, hence hyperbolicity of the quotient
$G_0/\langle R\rangle$. Here for growth we are rather interested in
the number of different elements of $G_0$ that can appear as boundary
words of fulfillable a abstract diagrams with given boundary
length (remember that our goal is to evaluate the number of equalities
$x=y$ holding in $G$ but not in $G_0$, with $x$ and $y$ elements of norm
at most $L$).

\subsection{Evaluation of growth}

We now turn back to random quotients: $R$ is a set of $(2m-1)^{d\ell}$
randomly chosen elements of $B_\ell$. Recall we set $L=A'\ell/2$ for some
value of $A'$ ensuring that if we know that $\abs{\B_L}\geq
(2m-1)^{g(1-\eps/2)L}$ then we know that the growth exponent of
$G=G_0/\langle R\rangle$ is at least $g(1-\eps)$.

We want to get an upper bound for the number $N$ of couples $x,y\in B_L$
such that $x=y$ in $G$ but $x\neq y$ in $G_0$. For any such couple there
is a worth-considering van Kampen diagram $D$ with boundary word
$xy^{-1}$, of boundary length at most $A'\ell$, with at most $K=A'/C$ new
relators, and at least one new relator (otherwise the equality $x=y$
would already occur in $G_0$). Let $\D$ be the decorated abstract van
Kampen diagram associated to $D$. Note that we have to count the number
of different couples $x,y\in B_L$ and \emph{not} the number of different
boundary words of van Kampen diagrams: since each $x$ and $y$ may have
numerous different representations as a word, the latter is higher than the
former.

We will show that, with overwhelming probability, we have $N\leq
\frac12(2m-1)^{gL}$.

The up to now free parameters $\eps_1$ and $\eps_2$ (in the definitions
of decorated abstract van Kampen diagrams and of apparent length) will be
fixed in the course of the proof, depending on $G_0$, $g$ and $d$ but not
on $\ell$. The length $\ell$ upon which our argument works will be set
depending on everything including $\eps_1$ and $\eps_2$.

\paragraph{Further notations.}
Let $n$ be the number of distinct relators in $\D$. For $1\leq a\leq n$,
let $m_a$ be the number of times the $a$-th relator appears in $\D$. Up
to reordering, we can suppose that the $m_a$ 's are non-increasing. Also
to avoid trivialities take $n$ minimal so that $m_n\geq 1$.

Let also $P_a$ be the probability that, if $a$ words $r_1,\ldots,r_a$ of
length $\ell$ are picked at random, there exist $n-a$ words
$r_{a+1},\ldots,r_n$ of lengt $\ell$ such that the relators
$r_1,\ldots,r_n$ fulfill $\D$. The $P_a$ 's are of course a
non-increasing sequence of probabilities. In particular, $P_n$ is the
probability that a random $n$-tuple of relators fulfills $\D$.

Back to our set $R$ of $(2m-1)^{d\ell}$ randomly chosen relators. Let $P^a$
be the probability that there exist $a$ relators $r_1,\ldots,r_a$ in
$R$, such that there exist words $r_{a+1}, \ldots,r_n$ of length $\ell$
such that the relators $r_1,\ldots,r_n$ fulfill $\D$. Again the $P^a$ 's
are a non-increasing sequence of probabilities and of course we have
\[
P^a\leq (2m-1)^{ad\ell} P_a
\]
since the $(2m-1)^{ad\ell}$ factor accounts for the choice of the $a$-tuple
of relators in $R$.

The probability that there exists a van Kampen diagram $D$ with respect
to the random presentation $R$, such that $\D$ is associated to $D$, is
by definition less than $P^a$ for any $a$. In particular, if for some
$\D$ we have $P^a\leq (2m-1)^{-\eps'\ell}$, then with overwhelming
probability,
$\D$ is not
associated to any van Kampen diagram of the random presentation. Since,
by Remark~\ref{numberofdavKds}, the number of possibilities for $\D$
grows subexponentially with $\ell$, we can sum this over $\D$ and
conclude that for any $\eps'>0$, with overwhelming probability
(depending on $\eps'$), all decorated abstract van Kampen diagrams $\D$ associated to
some van Kampen diagram of the random presentation satisfy $P^a\geq
(2m-1)^{-\eps'\ell}$ and in particular
\[
P_a \geq (2m-1)^{-ad\ell-\eps'\ell}
\]
which we assume from now on.

\medskip

We need to define one further quantity. Keep the notations of
Definition~\ref{defdavKds}. Let $1\leq a \leq n$ and let $1\leq i \leq k$
where $k$ is the number of translators of $\D$. Say that the $i$-th
translator is half finished at time $a$ if $r^\D(j_{2i})\leq a$ and
$r^\D(j_{2i+1})>a$, that is, if one side of the translator is a subword
of a relator $r_{a'}$ with $a'\leq a$ and the other of $r_{a''}$ with
$a''>a$. Now let $A_a$ be the sum of the apparent lengths of all
translators which are half finished at time $a$. In particular, $A_n$ is
the sum of the apparent lengths of all subwords $2i$ such that $2i$ is an
internal subword and $2i+1$ is a boundary subword of $\D$.

\paragraph{The proof.}
In this context, equation $(\star)$ (section~6.8) of~\cite{Oll03} reads
\[
A_a-A_{a-1}\geq
m_a\left(\ell(1-\eps'')+\frac{\log_{2m-1}P_a-\log_{2m-1}P_{a-1}}{\beta}\right)
\]
where $\eps''$ tends to $0$ when our free parameters $\eps_1,\eps_2$ tend
to $0$ (and $\eps''$ also absorbs the $o(\ell)$ term in~\cite{Oll03}). Also
recall that in the model of random quotient by random elements of balls we
have
\[
\beta=g/2
\]
by Proposition~20 of~\cite{Oll03}.

Setting $d'_a=\log_{2m-1}P_a$ and summing over $a$ we get, using $\sum
m_a=\abs{\D}$, that
\begin{eqnarray*}
A_n&\geq&
\left(\sum m_a\right)\ell \left(1-\eps''\right)+\frac2g \sum m_a(d'_a-d'_{a-1})
\\&=&\abs{\D}\ell(1-\eps'') +\frac2g \sum d'_a(m_a-m_{a+1})
\end{eqnarray*}

Now recall we saw above that for any $\eps'>0$, taking $\ell$ large enough
we can suppose that $P_a \geq (2m-1)^{-ad\ell-\eps'\ell}$, that is,
$d'_a+ad\ell+\eps'\ell\geq 0$. Hence
\begin{eqnarray*}
A_n&\geq&
\abs{\D}\ell(1-\eps'') + \frac2g \sum
(d'_a+ad\ell+\eps'\ell)(m_a-m_{a+1})
\\& &-\frac2g \sum
(ad\ell+\eps'\ell)(m_a-m_{a+1})
\\&=&
\abs{\D}\ell(1-\eps'') +\frac2g
\sum(d'_a+ad\ell+\eps'\ell)(m_a-m_{a+1})-\frac{d\ell}{g/2}\sum
m_a-\frac{\eps'\ell}{g/2} m_1
\\&\geq& \abs{\D}\ell(1-\eps'') +\frac{d'_n+nd\ell+\eps'\ell}{g/2}m_n
-\frac{d\ell+\eps'\ell}{g/2}\sum m_a
\end{eqnarray*}
where the last inequality follows from the fact that we chose the order
of the relators so that $m_a-m_{a+1}\geq 0$.

So using $m_n\geq 1$ we finally get
\[
A_n\geq
\abs{\D}\ell\left(1-\eps''-\frac{d+\eps'}{g/2}\right)+\frac{d'_n+nd\ell}{g/2}
\]

Set $\alpha=g/2-d>0$ so that this rewrites
\[
A_n\geq \frac2g\left(\abs{\D}\ell\left(\alpha-\eps'-\eps''g/2\right)
+d'_n+nd\ell
\right)
\]

Suppose the free parameters $\eps_1$, $\eps_2$ and $\eps'$ are chosen
small enough so that $\eps'+\eps''g/2\leq \alpha/2$ (remember that
$\eps''$ is a function of $\eps_1,\eps_2$ and $K$, tending to $0$ when
$\eps_1$ and $\eps_2$ tend to $0$). Since $\abs{\D}\geq 1$ (because we
are counting diagrams expressing equalities not holding in $G_0$) we get
$A_n\geq \ell\alpha/g+\frac2g \left(d'_n+nd\ell\right)$.

\begin{prop}
\label{alntuplesineq}
With overwhelming probability,
we can suppose that any
decorated abstract van Kampen diagram $\D$ satisfies
\[
A_n(\D)\geq\ell\alpha/g+\frac2g \left(d'_n(\D)+nd\ell\right)
\]
where $\alpha=g/2-d>0$.
\end{prop}

Let us translate back this inequality into a control on the numbers of
$n$-tuples of relators fulfilling $\D$.

Remember that, by definition, $d'_n$ is the log-probability that $n$
random relators $r_1,\ldots,r_n$ fulfill $\D$. As there are
$(2m-1)^{nd\ell}$ $n$-tuples of random relators in $R$ (by definition of
the density model), by linearity of expectation the expected number of
$n$-tuples of relators in $R$ fulfilling $\D$ is $(2m-1)^{nd\ell+d'_n}$.

By the Markov inequality, for given $\D$ the probability to pick a random
set $R$ such that the number of $n$-tuples of relators of $R$ fulfilling
$\D$ is greater than $(2m-1)^{nd\ell+d'_n+\eps'\ell}$, is less than
$(2m-1)^{-\eps'\ell}$. By Remark~\ref{numberofdavKds} the number of
possibilities for $\D$ is subexponential in $\ell$, and so, using
Proposition~\ref{alntuplesineq} we get

\begin{prop}
\label{numberntuples}
With overwhelming probability, we can suppose that for
any decorated abstract van Kampen diagram $\D$, the number of $n$-tuples
of relators in $R$ fulfilling $\D$ is at most 
\[
(2m-1)^{-\alpha\ell/2+g A_n(\D)/2+\eps'\ell}
\]
\end{prop}

Let us now turn back to the evaluation of the number of elements $x,y$ in
$B_L\subset G_0$ forming a van Kampen diagram $D$ with boundary word $xy^{-1}$.
For each such couple $x,y$ fix some geodesic writing of $x$ and $y$ as
words.
We will
first suppose that the abstract diagram $\D$ associated to $D$ is fixed
and evaluate the number of possible couples $x,y$ in function of $\D$, and then,
sum over the possible abstract diagrams $\D$.

So suppose $\D$ is fixed. Recall Proposition~\ref{propdavKd}: the
boundary word of $D$ is determined by giving two words for each
boundary-boundary translator, and one word for each internal-boundary
translator, this last one being subject to the apparent length condition
imposed in the definition of $\D$. By Remark~\ref{orderboundary}, the
number of ways to combine these subwords into a boundary word for $D$ is
controlled by $K$ and $\eps_2$ (independently of $\ell$).

In all the sequel, in order to avoid heavy notations, the notation
$\eps^\star$ will denote some function of $\eps'$, $\eps_1$ and $\eps_2$,
varying from time to time, and increasing when needed. The important
point is that $\eps^\star$ tends to $0$ when $\eps'$, $\eps_1$, $\eps_2$
do.

\bigskip

Let $(x_{2i},x_{2i+1})$ be a translator in $D$.  The definition of
translators
implies that there exist short words $\delta_1,\delta_2$, of length at
most $\eps_2(\abs{x_{2i}}+\abs{x_{2i+1}})$, such that
$x_{2i}\delta_1x_{2i+1}\delta_2=e$ in $G_0$. The words $x_{2i}$ and
$x_{2i+1}$ are either subwords of the geodesic words $x$ and $y$ making
the boundary of $D$, or subwords of relators in $R$; by definition of the
geodesic model of random quotient, the relators are geodesic as well. So
in either case $x_{2i}$ and $x_{2i+1}$ are geodesic\footnote{Except maybe
in the case when the translator straddles the end of $x$ and the
beginning of $y$ or conversely, or when it straddles the beginning and
end of a relator; these cases can be treated immediately by further
subdividing the translator, so we ignore this problem.}. Thus, the equality
$x_{2i}\delta_1x_{2i+1}\delta_2=e$ implies that $\norm{x_{2i+1}}\leq
\norm{x_{2i}}(1+\eps^\star)$ and conversely.
Also, by Remark~\ref{nosmalltrans}, we can suppose that
$\norm{x_{2i}}+\norm{x_{2i+1}}\geq \ell\eps_1\eps_2/6$, hence
$\norm{x_{2i}}\geq \ell\eps_1\eps_2(1-\eps^\star)/12$.

By definition of the growth exponent, there is some length $\ell_0$
depending only on $G_0$ such that if $\ell'_0\geq \ell_0$, then the
cardinal of $B_{\ell'_0}$ is at most $(2m-1)^{g(1+\eps')\ell'_0}$. So, if
$\ell$ is large enough (depending on $G_0$, $\eps_1$, $\eps_2$ and
$\eps'$) to ensure that $\ell\eps_1\eps_2(1-\eps^\star)/12\geq \ell_0$, we
can apply such an estimate to any $x_{2i}$.

\bigskip

To determine the number of possible couples $x,y$, we have to determine
the number of possibilites for each boundary-boundary or
internal-boundary translator $(x_{2i},x_{2i+1})$ (since by definition
internal translators do not contribute to the boundary).

First suppose that $(x_{2i},x_{2i+1})$ is a boundary-boundary translator.
Knowing the constraint $x_{2i}\delta_1x_{2i+1}\delta_2=e$, if $x_{2i}$
and $\delta_{1,2}$ are given then $x_{2i+1}$ is determined (as an element
of $G_0$). The number of
possibilities for $\delta_1$ and $\delta_2$ is at most
$(2m-1)^{2\eps_2(\norm{x_{2i}}+\norm{x_{2i+1}})}$. The number of
possibilities for $x_{2i}$ is at most $(2m-1)^{g(1+\eps')\norm{x_{2i}}}$
which, since $\norm{x_{2i}}\leq \frac12
\left(\norm{x_{2i}}+\norm{x_{2i+1}}\right) (1+\eps^\star)$, is
at most $(2m-1)^{\frac{g}{2}\left(\norm{x_{2i}}+\norm{x_{2i+1}}\right)
(1+\eps^\star)}$. So the total number of possibilities
for a boundary-boundary translator $(x_{2i},x_{2i+1})$ is at most
\[
(2m-1)^{\frac{g}{2}\left(\norm{x_{2i}}+\norm{x_{2i+1}}\right)(1+\eps^\star)}
\]
where of course the feature to remember is that the exponent is basically
$g/2$ times the total length $\norm{x_{2i}}+\norm{x_{2i+1}}$ of the translator.

Now suppose that $(x_{2i},x_{2i+1})$ is an internal-boundary translator.
The word $x_{2i}$ is by definition a subword of some relator $r_i\in
R$. So if a set of relators fulfilling $\D$ is fixed then $x_{2i}$ is
determined (we will multiply later by the number of possibilities for the
relators, using Proposition~\ref{numberntuples}). As above, the number of
possibilities for $\delta_1$ and $\delta_2$ is at most
$(2m-1)^{\eps^\star \norm{x_{2i}}}$. Once $x_{2i}$,
$\delta_1$ and $\delta_2$ are given, then $x_{2i+1}$ is determined (as an
element of $G_0$). So, if a set of relators fulfilling $\D$ is fixed,
then the number of possibilities for $x_{2i+1}$ is at most
$(2m-1)^{\eps^\star\norm{x_{2i}}}$, which reflects the fact that the set
of relators essentially determines the internal-boundary translators.

Let $A'_n$ be the sum of $\norm{x_{2i+1}}$ for all internal-boundary
translators $(x_{2i},x_{2i+1})$.  Let $B$ be the sum of
$\norm{x_{2i}}+\norm{x_{2i+1}}$ for all boundary-boundary translators. By
definition we have $\abs{\d \D}=A'_n+B$ maybe up to $\eps_1K\ell$.

So if a set of relators fulfilling $\D$ is fixed, then the total number
of possibilities for the boundary of $D$ is at most
\[
(2m-1)^{\frac{g}{2}\,B\,(1+\eps^\star)+\eps^\star A'_n}
\]
which, since both $B$ and $A'_n$ are at most $K\ell$, is at most
\[
(2m-1)^{gB/2+K\ell\eps^\star}
\]
(note that $A'_n$ does not come into play, since once the relators
fulfilling $\D$ are given, the internal-boundary translators are essentially determined).

The number of possibilities for an $n$-tuple of relators fulfilling $\D$
is given by Proposition~\ref{numberntuples}: it is at most
$(2m-1)^{-\alpha\ell/2+gA_n/2+\eps^\star\ell}$ (remember $\alpha=g/2-d$),
so that the total number of possibilities for the boundary of $D$ is at
most
\[
(2m-1)^{-\alpha\ell/2+(B+A_n)g/2+K\ell\eps^\star}
\]

Remember that $A_n$ is the sum of $\L(x_{2i})$ for all internal-boundary
translators $(x_{2i},x_{2i+1})$. By definition of apparent length we have
$\L(x_{2i})\leq \norm{x_{2i}}$.
Since in an internal-boundary translator $(x_{2i},x_{2i+1})$ we have
$\norm{x_{2i}}\leq \norm{x_{2i+1}}(1+\eps^\star)$, we get, after summing
on all internal-boundary translators, that $A_n\leq A'_n+K\ell\eps^\star$.
In particular, the above is at most
\[
(2m-1)^{-\alpha\ell/2+(B+A'_n)g/2+K\ell\eps^\star}
\]

Now remember that by definition we have $\abs{\d \D}=B+A'_n$ maybe up to
$\eps_1K\ell$ so that the above is in turn at most
\[
(2m-1)^{-\alpha\ell/2+\abs{\d \D}g/2+K\eps^\star\ell}
\]

This was for one decorated abstract van Kampen diagram $\D$. But by
Remark~\ref{numberofdavKds}, the number of such diagrams is
subexponential in $\ell$ (for fixed $K$ and $\eps_2$), and so, up to
increasing $\eps^\star$, this estimate holds for all diagrams
simultaneously.

\subsection{Conclusion}

Remember the discussion in the beginning of section~\ref{mainsection}. We
wanted to show that the cardinal $\abs{\B_L}$ of the ball of radius $L$
in $G$ was at least $(2m-1)^{gL(1-\eps/2)}$ for some $\eps$ chosen at the
beginning of our work.

We just proved that the number $N$ of couples of elements $x,y$ in $B_L$ such that there exists a van
Kampen diagram expressing the equality $x=y$ in $G$, but such that $x\neq
y$ in $G_0$ (which was expressed in the above argument by using that $D$
had at least one new relator) is at most
\[
(2m-1)^{-\alpha\ell/2+(\norm{x}+\norm{y})g/2+K\eps^\star\ell}
\]
where $\alpha=g/2-d>0$.

Now fix the free parameters $\eps'$, $\eps_1$, $\eps_2$ so that
$K\eps^\star\leq \alpha/4$ (this depends on $K$ and $G_0$ but not on
$\ell$; $K$ itself depends only on $G_0$). Choose $\ell$ large enough so
that all the estimates used above (implying every other variable) hold.
Also choose $\ell$ large enough (depending on $d$) so that $(2m-1)^{-\alpha\ell/4}\leq 1/2$.
We get
\[
N\leq \frac12 (2m-1)^{(\norm{x}+\norm{y})g/2}\leq \frac12 (2m-1)^{gL}
\]
since by assumption $\norm{x}$ and $\norm{y}$ are at most $L$. But on the
other hand we have $\abs{B_L}\geq (2m-1)^{gL}$ and so
\[
\abs{\B_L}\geq \abs{B_L}-N \geq \frac12 (2m-1)^{gL} \geq (2m-1)^{gL(1-\eps/2)}
\]
as soon as $\ell$ is large enough (since $L$ is a multiple of $\ell$), which ends the proof.

\appendix
\renewcommand{\thesection}{Appendix:}

\section{Locality of growth in hyperbolic groups}

The goal of this section is to show that, in a hyperbolic group, if we
know an estimate of the growth exponent in some finite ball of the group,
then this provides an estimate of the growth exponent of the group (whose
quality depends on the radius of the given finite ball).

Let $G=\presgroup{a_1,\ldots,a_m}{R}$ be a $\delta$-hyperbolic group
generated by the elements $a_i^{\pm 1}$, with $m\geq 2$. For $x\in G$ let
$\norm{x}$ be the norm of $x$ with respect to this generating set. Let
$B_\ell$ be the set of elements of norm at most $\ell$.

\begin{prop}
\label{localgr}

Suppose that for some $g>0$, for some $\ell_0\geq 2\delta+4/g$ and
$\ell_1\geq A\ell_0$, with $A\geq 500$, we
have
\[
\abs{B_{\ell_0}}\leq (2m-1)^{1.1g\ell_0}
\]
and
\[
\abs{B_{\ell_1}}\geq (2m-1)^{g\ell_1}
\]

Then the growth exponent of $G$ is at least $g(1-40/A)$.
\end{prop}

Note that the occurrence of $1/g$ in the scale upon which the proposition
is true is natural: indeed, an assumption such as $\abs{B_\ell}\geq
(2m-1)^{g\ell}$ for $\ell<1/g$ is not very strong... The growth $g$ can be
thought of as the inverse of a length, so this result is homogeneous.

\begin{cor}
The growth exponent of a presentation of a hyperbolic group is
computable.
\end{cor}

\begin{dem}
Indeed, remember from~\cite{Pap96} (after~\cite{Gro87}) that the
hyperbolicity constant $\delta$ of a presentation of a hyperbolic group
is computable. Thanks to the isoperimetric
inequality, the word problem in a hyperbolic group is solvable, so that for
any $\ell$ an exact computation of the cardinal of $B_\ell$ is possible.
Setting $g_\ell=\frac1\ell \log_{2m-1}\abs{B_\ell}$, we know that $g_\ell$
will converge to some (unknown) positive value, so that 
$g_{\ell}$ and $g_{A\ell}$ will become arbitrarily close, and since
$g_\ell$ is bounded from below sooner or later we will have $\ell\geq
2\delta+4/g_{A\ell}$, in which case we can apply the proposition to
$\ell$ and $A\ell$.
\end{dem}


\begin{dem}[ of the proposition]

Let $(,)$ denote the Gromov product in $G$, with origin at $e$, that is
\[
(x,y)=\frac12 \left(\norm{x}+\norm{y}-\norm{x-y}\right)
\]
for $x,y\in G$, where, following~\cite{GH}, we write $\norm{x-y}$ for
$\norm{x^{-1}y}=\norm{y^{-1}x}$. Since triangles are $\delta$-thin, we
have (\cite{GH}, Proposition~2.21) for any three points $x$,$y$, $z$ in
$G$
\[
(x,z)\geq \min\left((x,y),(y,z)\right)-2\delta
\]

Let $S_\ell$ denote the set of elements of norm $\ell$ in the hyperbolic
group $G$.  Consider also, for homogeneity reasons, the annulus
$S_{\ell,a}=B_\ell\setminus B_{\ell-a}$.


\begin{prop}
Let $g\in B_\ell$ and let $a\geq 0$. The number of elements $g'$ in $S_\ell$
or $B_\ell$ such that $(g,g')\geq a$ is at most
$\abs{B_{\ell-a+2\delta}}$.
\end{prop}

\begin{dem}
Suppose that $(g,g')\geq a$.
Let $x$ be the point at distance $a$ from $e$ on some geodesic joining
$e$ to $g$. By construction we have $(g,x)=a$. But
\[
(g',x)\geq \min\left((g',g),(g,x)\right)-2\delta \geq a-2\delta
\]
and unwinding the definition of $(g',x)$ yields
\[
\norm{g'-x}\leq \norm{g'}+\norm{x}-2a+2\delta\leq\ell-a+2\delta
\]

So $g'$ lies at distance at most $\ell-a+2\delta$ from $x$, hence the
number of possibilities for $g'$ is at most $\abs{B_{\ell-a+2\delta}}$.
(This is most clear on a picture.)
\end{dem}

We know show that, if we multiply two elements of the sphere $S_\ell$
then we often get an element of norm close to $2\ell$.

\begin{cor}\label{multgeod}
Let $g\in S_{\ell,a}$. The number of elements $g'$ in $S_{\ell,a}$ such that
$\norm{gg'}\geq 2\ell-4a$ is at least
$\abs{S_{\ell,a}}-\abs{B_{\ell-a+2\delta}}$.
\end{cor}

\begin{dem}
We have $\norm{gg'}=\norm{g}+\norm{g'}-2(g^{-1},g')$. So if $\norm{g}\geq
\ell-a$, $\norm{g'}\geq \ell-a$ and $(g^{-1},g')\leq a$,
then $\norm{gg'}\geq 2\ell-4a$.

But by the last proposition, the number of ``bad'' elements $g'$ such that
$(g^{-1},g')\geq a$ is at most $\abs{B_{\ell-a+2\delta}}$.
\end{dem}

So multiplying long elements often gives twice as long elements. We now
show that this procedure does not build too often the same new element.

\begin{prop}
Let $x\in S_{2\ell, 4a}$. The number of couples $(g,g')$ in
$S_{\ell,a}\times S_{\ell,a}$ such that $x=gg'$ is at most
$\abs{B_{6a+2\delta}}$.
\end{prop}

\begin{dem}
Choose a geodesic decomposition $x=hh'$ with
$\norm{h}=\norm{h'}=\norm{x}/2$. It is easy to see that if $x=gg'$ as
above, then $g$ is $6a+2\delta$-close to $h$ (and then $g'$ is
determined).
\end{dem}

Combining the last two results yields the following ``almost
supermultiplicative'' estimate for the cardinals of balls (compare the
trivial converse inequality $\abs{B_{2\ell}}\leq \abs{B_\ell}^2$).

\begin{cor}\label{growthmult}
\[
\abs{B_{2\ell}}\geq
\frac{1}{\abs{B_{6a+2\delta}}}\left(\abs{B_\ell}-2\abs{B_{\ell-a+2\delta}}\right)^2
\]
\end{cor}

\begin{dem}
Indeed, the last two results imply that
\[
\abs{S_{2\ell,4a}}\geq \frac{1}{\abs{B_{6a+2\delta}}}
\abs{S_{\ell,a}}\left(\abs{S_{\ell,a}}-\abs{B_{\ell-a+2\delta}}\right)
\]
which implies the above by the trivial estimates $\abs{B_{2\ell}}\geq
\abs{S_{2\ell,4a}}$ and $\abs{S_{\ell,a}}\geq
\abs{B_\ell}-\abs{B_{\ell-a+2\delta}}$.
\end{dem}

In order to apply this, we need to know both that $\abs{B_\ell}$ is big
and that $\abs{B_{\ell-a}}$ is not too big compared to $\abs{B_\ell}$.
Asymptotically one would expect $\abs{B_{\ell-a}}\approx
(2m-1)^{-ga}\abs{B_\ell}$. The next lemma states that, under the
assumptions of Proposition~\ref{localgr}, we can almost
realize this, up to changing $\ell$ by some controlled factor.

\begin{lem}\label{lemcrown}
Suppose that for some $g$, for some $\ell_0$ and $\ell_1\geq 100\ell_0$
we have $\abs{B_{\ell_0}}\leq (2m-1)^{1.2g\ell_0}$ and
$\abs{B_{\ell_1}}\geq (2m-1)^{g\ell_1}$.
Let $a\leq\ell_0$. There exists $0.65\ell_1\leq \ell\leq \ell_1$ such that
\[
\abs{B_\ell}\geq (2m-1)^{g\ell}
\]
and
\[
\abs{B_\ell}\geq (2m-1)^{ga/2}\abs{B_{\ell-a}}
\]
\end{lem}

\begin{dem}[ of the lemma]
First, note that by subadditivity, the inequality $\abs{B_{\ell_0}}\leq
(2m-1)^{1.2g\ell_0}$ implies that for any $\ell$, writing $\ell=k\ell_0-r$
($k\in \N, 0\leq r<\ell_0$) we have $\abs{B_{\ell}}\leq
(2m-1)^{1.2kg\ell_0}$. Especially for $\ell\geq 50\ell_0$ we have
$1\leq k\ell_0/\ell\leq 51/50$ and so in particular, if
$\ell_1\geq 100\ell_0$ then $\abs{B_{0.65\ell_1}}\leq (2m-1)^{0.8g\ell_1}$
(indeed $0.65\times1.2\times 51/50\leq 0.8$).

Suppose that for all $0.65\ell_1\leq \ell\leq \ell_1$ with $\ell=\ell_1-ka$
($k\in \N$) we have
$\abs{B_\ell}< (2m-1)^{ga/2}\abs{B_{\ell-a}}$. Write
$\ell_1-0.65\ell_1=qa-r$ with $q\in\N$, $0\leq r<a$.
Then we get
\begin{eqnarray*}
\abs{B_{\ell_1}}
&<&(2m-1)^{ga/2}\abs{B_{\ell_1-a}}< (2m-1)^{ga}\abs{B_{\ell_1-2a}}< \cdots
\\&<&
(2m-1)^{gqa/2}\abs{B_{0.65\ell_1-r}}\leq
(2m-1)^{g(\ell_1-0.65\ell_1)/2+ga/2}\abs{B_{0.65\ell_1}}
\\&\leq&
(2m-1)^{g(0.35\ell_1)/2+g\ell_1/200+0.8g\ell_1} < (2m-1)^{0.98g\ell_1}
\end{eqnarray*}
contradicting the assumption.

So we can safely take the largest $\ell\leq \ell_1$ satisfying $\abs{B_\ell}\geq
(2m-1)^{ga/2}\abs{B_{\ell-a}}$ and such that $\ell_1-\ell$ is a multiple of
$a$.

Since $\ell$ is largest, for $\ell\leq \ell'\leq \ell_1$
we have $\abs{B_{\ell'}}\leq (2m-1)^{ga/2}\abs{B_{\ell'-a}}$. We get,
$a$-step by $a$-step, that
$\abs{B_{\ell_1}}\leq (2m-1)^{g(\ell_1-\ell)/2}\abs{B_\ell}$. Using
the assumption $\abs{B_{\ell_1}}\geq (2m-1)^{g\ell_1}$ we now get
$\abs{B_\ell}\geq(2m-1)^{g\ell_1-g(\ell_1-\ell)/2}\geq (2m-1)^{g\ell}$ as
needed.
\end{dem}

Now equipped with the lemma, we can apply Corollary~\ref{growthmult} to
show that if we know that $B_\ell$ is large for some $\ell$, then we get
a larger $\ell'$ such that $B_{\ell'}$ is large as well. We will then
conclude by induction.

\begin{lem}
Suppose that for some $g$, for some $\ell_0\geq 2\delta+4/g$ and
$\ell_1\geq A\ell_0$ (with $A\geq 100$)
we have $\abs{B_{\ell_0}}\leq (2m-1)^{1.2g\ell_0}$ and
$\abs{B_{\ell_1}}\geq (2m-1)^{g\ell_1}$.
Then there exists $\ell_2\geq 1.3\ell_1$ such that
\[
\abs{B_{\ell_2}}\geq (2m-1)^{g\ell_2(1-9/A)}
\]
\end{lem}

\begin{dem}[ of the lemma]
Consider the $\ell$ provided by Lemma~\ref{lemcrown} where we take
$a=\ell_0$. This provides an $\ell\geq 0.65\ell_0$ such that
$\abs{B_\ell}\geq (2m-1)^{g\ell}$
and $\abs{B_\ell}\geq (2m-1)^{ga/2}\abs{B_{\ell-a}}$.

So by Corollary~\ref{growthmult} (applied to $2a$ instead of $a$) we have
\[
\abs{B_{2\ell}}\geq \frac{1}{\abs{B_{12a+2\delta}}} \abs{B_\ell}^2
\left(1-2\abs{B_{\ell-2a+2\delta}}/\abs{B_\ell}\right)^2
\]
Since $a=\ell_0\geq 2\delta$ we have $\ell-2a+2\delta\leq \ell-\ell_0$
and so
\[
\abs{B_{2\ell}}\geq \frac{1}{\abs{B_{12\ell_0+2\delta}}} \abs{B_\ell}^2
\left(1-2(2m-1)^{-g\ell_0/2}\right)^2
\]

If $\ell_0\geq 4/g$, since $2m-1\geq 2$ we have
$\left(1-2(2m-1)^{-g\ell_0/2}\right)^2\geq 1/4$ and so
\[
\abs{B_{2\ell}}\geq\frac{1}{4\abs{B_{12\ell_0+2\delta}}}\abs{B_\ell}^2
\]

We have $\abs{B_{12\ell_0+2\delta}}\leq \abs{B_{13\ell_0}}\leq
\abs{B_{\ell_0}}^{13}$ by subadditivity. So by the assumptions
\[
\abs{B_{2\ell}}\geq \frac{1}{4\abs{B_{\ell_0}}^{13}}\abs{B_\ell}^2
\geq (2m-1)^{2g\ell-16g\ell_0-2}=(2m-1)^{2g\ell(1-8\ell_0/\ell-1/g\ell)}
\]
which is at least $(2m-1)^{2g\ell(1-9/A)}$ since $8\ell_0/\ell\leq 8/A$
and $1/g\ell\leq1/gA\ell_0\leq 1/A$ since $\ell_0\geq 4/g$.

So we can take $\ell_2=2\ell$, which is at least $1.3\ell_1$.
\end{dem}

Now the proposition is clear: start from $\ell_1$ and construct by
induction a sequence $\ell_i$ with $\ell_{i+1}\geq 1.3 \ell_i$ using the
lemma applied to $\ell_0$ and $\ell_i$; thus
\[
\abs{B_{\ell_i}}\geq (2m-1)^{g\ell_i\prod_{k=0}^{i-2} (1-9/(A\cdot1.3^k))}
\]
and note that the infinite product converges to a value greater than
$1-40/A$.
The only thing to check is that, in order to be allowed to
apply the previous lemma to $\ell_0$ and $\ell_i$ at each step, we must ensure that
$1.1/(1-40/A)\leq 1.2$, which is guaranteed as soon as $A\geq 500$.
\end{dem}



\newpage

\begin{thebibliography}{MMM}

\bibitem[Ch93]{Ch93} C.~Champetier, \emph{Cocroissance des groupes à petite
simplification}, Bull.\ London Math.\ Soc.~\textbf{25} (1993), No.~5,
438--444.

\bibitem[Ch95]{Ch95} C.~Champetier, \emph{Propriétés statistiques des
groupes de présentation finie}, J.\ Adv.\ Math.\ \textbf{116} (1995),
No.~2, 197--262.

\bibitem[Gh]{Gh} É.~Ghys, \emph{Groupes aléatoires}, séminaire
Bourbaki~\textbf{916} (2003).

\bibitem[GhH90]{GH} É.~Ghys, P.~de la Harpe, \emph{Sur les groupes
hyperboliques d'après Mikhael Gromov}, Progress in Math.~\textbf{83},
Birkhäuser (1990).

\bibitem[GrH97]{GdlH} R.I.~Grigorchuk, P.~de la Harpe, \emph{On problems
related to growth, entropy, and spectrum in group theory}, Dynam.\
Control Systems \textbf{3} (1997), No.~1, 51--89.

\bibitem[Gro87]{Gro87} M.~Gromov, \emph{Hyperbolic Groups}, in
\emph{Essays in group theory}, ed.\ S.M.~Gersten, Springer (1987),
75--265.

\bibitem[Gro93]{Gro93} M.~Gromov, \emph{Asymptotic Invariants of Infinite
Groups}, in \emph{Geometric group theory}, ed.\ G.~Niblo, M.~Roller,
Cambridge University Press, Cambridge (1993).

\bibitem[Gro03]{Gro03} M.~Gromov, \emph{Random Walk in Random Groups},
Geom.\ Funct.\ Anal.~\textbf{13} (2003), No.~1, 73--146.

\bibitem[Har00]{Har00} P.~de la Harpe, \emph{Topics in geometric group
theory}, Chicago University Press (2000).

\bibitem[Har02]{Har02} P.~de la Harpe, \emph{Uniform growth in groups of
exponential growth}, Geom.\ Dedicata~\textbf{95} (2002), 1--17.

\bibitem[HLS]{HLS} N.~Higson, V.~Lafforgue, G.~Skandalis,
\emph{Counterexamples to the Baum-Connes conjecture}, Geom.\ Funct.\
Anal.~\textbf{12} (2002), No.~2, 330--354.

\bibitem[Kou98]{Kou98} M.~Koubi, \emph{Croissance uniforme dans les groupes
hyperboliques}, Ann.\ Institut Fourier~\textbf{48} (1998), No.~5,
1441--1453.

\bibitem[LS77]{LS} R.C.~Lyndon, P.E.~Schupp, \emph{Combinatorial Group
Theory}, Ergebnisse der Mathematik und ihrer Grenzgebiete~\textbf{89},
Springer (1977).

\bibitem[Oll03]{Oll03} Y.~Ollivier, \emph{Sharp phase transition theorems for
hyperbolicity of random groups}, to appear in GAFA, Geom.\ Funct.\
Anal.

\bibitem[Oll04]{Oll04} Y.~Ollivier, \emph{Cogrowth and spectral gap of
generic groups}, submitted. Available on
ArXiv: \texttt{math.GR/0401048}

\bibitem[Ols92]{Ols92} A.Yu.~Ol'shanski\u{\i}, \emph{Almost Every Group is
Hyperbolic}, Int.\ J.\ Algebra Comput.~\textbf{2} (1992), No.~1, 1--17.

\bibitem[Pap96]{Pap96} P.~Papasoglu, \emph{An Algorithm Detecting
Hyperbolicity}, in G.~Baumslag (ed.) et al., \emph{Geometric and
Computational Perspectives on Infinite Groups}, DIMACS Ser.\ Discrete
Math.\ Theor.\ Comput.\ Sci.\ \textbf{25} (1996), 193--200.

\bibitem[Sho91]{Sho} H.~Short et al., in \emph{Group Theory from a Geometrical
Viewpoint}, ed.\ É.~Ghys, A.~Haefliger, A.~Verjovsky, World Scientific
(1991).

\bibitem[Shu99]{Shu} A.G.~Shukhov, \emph{On the dependence of the growth
exponent on the length of the defining relation}, Math.\
Notes~\textbf{65} (1999), No.~3--4, 510--515.

\bibitem[Ver00]{Ver00} A.M.~Vershik, \emph{Dynamic theory of growth in
groups: entropy, boundaries, examples}, Russian Math.\ Surveys~\textbf{55}
(2000), No.~4, 667--733.

\bibitem[Wil04]{Wil04} J.S.~Wilson, \emph{On exponential growth and
uniform exponential growth for groups}, Invent.\ Math~\textbf{155}
(2004), No.~2, 287--303.


\end{thebibliography}
\end{document}